\documentclass[letterpaper,12pt]{article}
\usepackage[utf8]{inputenc}
\usepackage[T1]{fontenc}
\usepackage{lmodern}
\usepackage[margin=1.0in]{geometry}
\usepackage[numbers,super]{natbib}
\usepackage{mathtools,amsmath,amssymb}
\usepackage{microtype}
\usepackage{appendix}
\usepackage[ruled]{algorithm2e}
\usepackage{hyperref}
\usepackage{setspace}      %singlespacing/doublespacing[onehalfspacing]
\allowdisplaybreaks
\usepackage{blindtext}

\title{Augmented Lagrangian method for a TV-based model for demodulating phase discontinuities}
\author{Ricardo Legarda-Saenz and Carlos Brito-Loeza\\
CLIR at Facultad de Matemáticas, Universidad Autónoma de Yucatán\\
Apartado Postal 172. 97110 Mérida, Yucatán. México\\
E-mail: \texttt{rlegarda@correo.uady.mx}  } 
\date{\today}

\begin{document}
\maketitle

\begin{abstract}
In this work we reformulate the method presented in App. Opt. 53:2297 (2014) as a constrained minimization problem using the augmented Lagrangian method. First we introduce the new method and then describe the numerical solution, which results in a simple algorithm. Numerical experiments with both synthetic and real fringe patterns show the accuracy and simplicity of the resulting algorithm.
\end{abstract}

\section{Introduction}
The main goal of fringe analysis techniques is to recover accurately the local modulated phase from one or several fringe patterns~\cite{Servin2014}; such phase is related to some physical quantities like shape, deformation, refractive index, temperature, etc. The basic model for a fringe pattern is given by
\begin{equation}\label{eq:Franjas} 
I_{\mathbf{x}} = a_{\mathbf{x}} + b_{\mathbf{x}}\cos\left(\omega_{\mathbf{x}} + \phi_{\mathbf{x}}\right)
\end{equation}
where $\mathbf{x} = (x_{1},x_{2})$, $a_{\mathbf{x}}$ is the background illumination, $b_{\mathbf{x}}$ is the amplitude modulation, and $\phi_{\mathbf{x}}$ is the phase map to be recovered; the spatial carrier frequency of the fringe pattern is defined by the term $\omega_{\mathbf{x}}$.

Several methods which successfully estimate the phase from a single pattern have been reported in literature.~\cite{Marroquin1999,Legarda-Saenz2002a,Villa2000,Rivera2005} These methods consider phase maps, amplitude and illumination terms as a continuous. However, the recovery of a discontinuous phase map from a single fringe pattern remains a pending task, and is a challenging problem. Some years ago, it was proposed a method for computing discontinuous phase maps of a fringe pattern with carrier frequency, based on the minimization of a regularized cost function, which uses a second-order edge-preserving potential~\cite{Galvan2006}. Although this method is able to detect and reconstruct phase discontinuities, its cost functional is not convex, hence convergence to an optimal solution is conditioned to the provided initial phase usually computed by standard methods.
 
In a recent work, it was proposed a method for computing discontinuous phase maps based on a total variational (TV) approach,~\cite{Rudin1992} where TV regularization is applied to the background, amplitude and phase terms of the fringe model, resulting in accurate phase reconstructions.~\cite{Legarda-Saenz2014} Despite this fact, this model lacks of a fast algorithm for its solution. Recently, it was proposed a fixed point method to speed up the numerical solution of this model.~\cite{Brito-Loeza2019} This fixed point method shows a good performance solving the model presented in reference \citenum{Legarda-Saenz2014}; however, similar to other methods based on the TV approach, the fixed point method performance is dramatically reduced for problems highly anosotropic.

In this work we reformulate the model presented in reference \citenum{Legarda-Saenz2014} as a constrained minimization problem using augmented Lagrangian method. First, we describe the ideas that give support to our reformulation,  and then we describe the numerical solution of the proposed augmented Lagrangian method, which results in a simple algorithm. The performance of the proposed method is evaluated by numerical experiments with both synthetic and real data. A comparison against the fixed point method is presented. Finally we discuss our results and present some concluding remarks.

\section{Methodology}

\subsection{Fixed point method for computing discontinuous phase maps based on TV model}\label{sec:fixedpoit}
The method proposed for computing discontinuous phase maps based on TV approach~\cite{Legarda-Saenz2014} is given by 
\begin{multline}\label{eq:energiaMinimiza}
\underset{\phi,a,b}{\min}\;E(\phi_{\mathbf{x}}, a_{\mathbf{x}}, b_{\mathbf{x}};\omega_{\mathbf{x}}) = \frac{\lambda}{2}\int_{\Omega}\left(I_{\mathbf{x}} - g_{\mathbf{x}}\right)^{2}d\mathbf{x} + \int_{\Omega}\left|\nabla\phi_{\mathbf{x}}\right| d\mathbf{x}\\ 
+ \int_{\Omega}\left|\nabla a_{\mathbf{x}}\right| d\mathbf{x} + \int_{\Omega}\left|\nabla b_{\mathbf{x}}\right| d\mathbf{x}
\end{multline}
where $\Omega\subset \mathbb{R}^{2}$ is the domain of integration, $g_{\mathbf{x}}$ is the given fringe pattern and $\lambda$ is the regularization parameter. As it was shown in ref. \citenum{Legarda-Saenz2014}, this model allows the accurate demodulation of a single fringe pattern with discontinuities. The first-order optimality conditions or Euler-Lagrange equations of Eq. (\ref{eq:energiaMinimiza}) are given by
\begin{align}\label{eq:gradienteE}
  \lambda\left(I_{\mathbf{x}} - g_{\mathbf{x}}\right)\,\frac{\partial I_{\mathbf{x}}}{\partial\phi_{\mathbf{x}}} 
  &-\nabla\cdot\frac{\nabla\phi_{\mathbf{x}}}{\left|\nabla\phi_{\mathbf{x}}\right|} = 0,\nonumber\\
  \lambda\left(I_{\mathbf{x}} - g_{\mathbf{x}}\right)\,\frac{\partial I_{\mathbf{x}}}{\partial a_{\mathbf{x}}}
  &- \nabla\cdot\frac{\nabla a_{\mathbf{x}}}{\left|\nabla a_{\mathbf{x}}\right|} = 0,\\
  \lambda\left(I_{\mathbf{x}} - g_{\mathbf{x}}\right)\,\frac{\partial I_{\mathbf{x}}}{\partial b_{\mathbf{x}}}
  &- \nabla\cdot\frac{\nabla b_{\mathbf{x}}}{\left|\nabla b_{\mathbf{x}}\right|} = 0,\nonumber
\end{align}
with boundary conditions
\begin{equation}
\frac{\nabla\phi_{\mathbf{x}}}{\left|\nabla\phi_{\mathbf{x}}\right|}\cdot\mathbf{n} = 0,\quad
\frac{\nabla a_{\mathbf{x}}}{\left|\nabla a_{\mathbf{x}}\right|}\cdot\mathbf{n} = 0,\quad
\frac{\nabla b_{\mathbf{x}}}{\left|\nabla b_{\mathbf{x}}\right|}\cdot\mathbf{n} = 0,
\end{equation}
where $\mathbf{n}$ denotes the unit outer normal to the boundary.

The numerical solution proposed in ref. \citenum{Legarda-Saenz2014} was a gradient descent scheme which is very slow and therefore a large number of iterations are necessary to reach an adequate solution. To speed up the convergence, a fixed point method~\cite{Vogel1996,Vogel2002} was proposed recently:~\cite{Brito-Loeza2019} 

First, Eq. (\ref{eq:gradienteE}) is written in the following way
\begin{align}\label{eq:gradienteE1} 
  \lambda\Big(a_{\mathbf{x}} + b_{\mathbf{x}}\cos\psi_{\mathbf{x}} - g_{\mathbf{x}}\Big)\Big(-b_{\mathbf{x}}\sin\psi_{\mathbf{x}}\Big) &-\nabla\cdot\frac{\nabla\phi_{\mathbf{x}}}{\left|\nabla\phi_{\mathbf{x}}\right|} = 0,\nonumber\\
  \lambda\Big(a_{\mathbf{x}} + b_{\mathbf{x}}\cos\psi_{\mathbf{x}} - g_{\mathbf{x}}\Big)\Big(\cos\psi_{\mathbf{x}}\Big) &- \nabla\cdot\frac{\nabla b_{\mathbf{x}}}{\left|\nabla b_{\mathbf{x}}\right|} = 0,\\  
  \lambda\Big(a_{\mathbf{x}} + b_{\mathbf{x}}\cos\psi_{\mathbf{x}} - g_{\mathbf{x}}\Big) &- \nabla\cdot\frac{\nabla a_{\mathbf{x}}}{\left|\nabla a_{\mathbf{x}}\right|} = 0,\nonumber
\end{align}
where the term $\psi_{\mathbf{x}}$ is defined as $\psi_{\mathbf{x}} = \omega_{\mathbf{x}} + \phi_{\mathbf{x}}.$ 

In ref. \citenum{Vogel1996} it was proposed a fixed point method to solve the TV model. The basic idea is to linearize the nonlinear term of this model, so that at each iteration $k$ the method is required to solve a linear system of the form \[L_{d}\left(d_{\mathbf{x}}^{k}\right)d_{\mathbf{x}}^{k+1} = f_{d},\quad k = 1, 2,\ldots\] where $d_{\mathbf{x}}$ is the unknown variable, the operator $L_{d}(\cdot)$ has been made linear by lagging the nonlinear term $\tfrac{1}{\left|\nabla d_{\mathbf{x}}\right|},$ and $ f_{d}$ has the terms which remain constant at each $k$ iteration. 

Arranging Eq. (\ref{eq:gradienteE1}) in the same way as described before, the proposed fixed point iteration is given by
\begin{align}\label{eq:fixedpoint1}
  -\lambda\left(b^{k}_{\mathbf{x}}\right)^{2}\sin\psi^{k}_{\mathbf{x}}\cos\psi^{k+1}_{\mathbf{x}} -\nabla\cdot\frac{\nabla\phi^{k+1}_{\mathbf{x}}}{\left|\nabla\phi^{k}_{\mathbf{x}}\right|} &= \lambda\Big(a^{k}_{\mathbf{x}} - g_{\mathbf{x}}\Big)\Big(b^{k}_{\mathbf{x}}\sin\psi^{k}_{\mathbf{x}}\Big),
  \nonumber\\  
  \Big(\lambda\cos^{2}\psi^{k}_{\mathbf{x}} - \nabla\cdot\frac{\nabla}{\left|\nabla b^{k}_{\mathbf{x}}\right|}\Big) b^{k+1}_{\mathbf{x}} &= -\lambda\Big(a^{k}_{\mathbf{x}} - g_{\mathbf{x}}\Big)\cos\psi^{k}_{\mathbf{x}},\\  
  \Big(\lambda - \nabla\cdot\frac{\nabla}{\left|\nabla a^{k}_{\mathbf{x}}\right|}\Big) a^{k+1}_{\mathbf{x}} &= -\lambda\Big(b^{k}_{\mathbf{x}}\cos\psi^{k}_{\mathbf{x}} - g_{\mathbf{x}}\Big),\nonumber
\end{align}
with $\psi^{k}_{\mathbf{x}} = \omega_{\mathbf{x}} + \phi^{k}_{\mathbf{x}}.$ 

As can be observed in the first term of Eq. (\ref{eq:fixedpoint1}), it is not possible to separate the term $\psi^{k+1}_{\mathbf{x}}$ from the cosine function. To fix this, the cosine function is linearised in the following way \[\cos\psi^{k+1}_{\mathbf{x}}\approx\cos\psi^{k}_{\mathbf{x}} - \left(\phi^{k+1}_{\mathbf{x}}-\phi^{k}_{\mathbf{x}}\right)\sin\psi^{k}_{\mathbf{x}}\]  and the proposed fixed point iteration is given by~\cite{Brito-Loeza2019} 
\begin{align}\label{eq:fixedpoint}
  &\Big(\lambda \left(b^{k}_{\mathbf{x}}\right)^{2}\sin^{2}\psi^{k}_{\mathbf{x}} - \nabla\cdot\frac{\nabla}{\left|\nabla\phi^{k}_{\mathbf{x}}\right|} \Big)\phi^{k+1}_{\mathbf{x}} = \lambda \Big(a^{k}_{\mathbf{x}}b_{\mathbf{x}}\sin\psi^{k}_{\mathbf{x}}\nonumber\\
  &+ \left(b^{k}_{\mathbf{x}}\right)^{2}\cos\psi^{k}_{\mathbf{x}}\sin\psi^{k}_{\mathbf{x}} + \phi^{k}_{\mathbf{x}}\left(b^{k}_{\mathbf{x}}\right)^{2}\sin^{2}\psi^{k}_{\mathbf{x}} - g_{\mathbf{x}}b_{\mathbf{x}}\sin\psi^{k}_{\mathbf{x}}\Big),\nonumber\\
  &\Big(\lambda\cos^{2}\psi^{k}_{\mathbf{x}} - \nabla\cdot\frac{\nabla}{\left|\nabla b^{k}_{\mathbf{x}}\right|}\Big) b^{k+1}_{\mathbf{x}} = -\lambda\Big(a^{k}_{\mathbf{x}}\cos\psi^{k}_{\mathbf{x}} - g_{\mathbf{x}}\cos\psi^{k}_{\mathbf{x}}\Big),\\  
  &\Big(\lambda - \nabla\cdot\frac{\nabla}{\left|\nabla a^{k}_{\mathbf{x}}\right|}\Big) a^{k+1}_{\mathbf{x}} = -\lambda\Big(b^{k}_{\mathbf{x}}\cos\psi^{k}_{\mathbf{x}} - g_{\mathbf{x}}\Big).\nonumber 
\end{align}
Ref. \citenum{Brito-Loeza2019} provides a detailed explanation of the convergence proof and the numerical performance of this method .

\subsection{Augmented Lagrangian method for TV model}\label{sec:almTV}
A distinctive feature of the solution proposed in Eq. (\ref{eq:fixedpoint}) is that all the PDE's have the coefficient $\tfrac{1}{\left|\nabla d_{\mathbf{x}}\right|},$ which is quite hard to deal with it numerically due to the inherent discontinuity. A typical solution to this problem is to include a small constant to avoid division by zero, that is $\tfrac{1}{\sqrt{\left|\nabla d_{\mathbf{x}}\right|^{2}+\beta}},$ but this affect both the accuracy and efficiency of the solution. Here, we will review an alternative to avoid this problem in a similar TV model. This  will prove to be helpful when we introduce our proposed augmented Lagrangian method in the next section. 

The TV model for image denoising is given by~\cite{Rudin1992}
\begin{equation}\label{eq:TVmethod}
\underset{u}{\min}\;F_{TV}(u) = \frac{\lambda}{2}\int_{\Omega}\left(u - f\right)^{2}d\mathbf{x} + \int_{\Omega}\left|\nabla u\right| d\mathbf{x},
\end{equation}
where $u$ is the original image and $f$ is the noisy image. It is well known that the computation of the TV method suffers from non-differentially due to the TV norm.~\cite{Rudin1992,Chan1999,Chambolle2004,Getreuer2012} Many numerical methods have been proposed to improve this drawback; one of them is to convert the above functional into a constrained optimization problem, where an auxiliary term is introduced to separate the estimation of the non-differentiable term:~\cite{Hestenes1969,Rockafellar1973,Bertsekas1996,Nocedal2006}
\begin{align}\label{eq:constrainedTV}
\underset{u, q}{\min}\;G_{TV}(u,q) &= \frac{\lambda}{2}\int_{\Omega}\left(u - f\right)^{2}d\mathbf{x} + \int_{\Omega}\left|q\right| d\mathbf{x}\\
&\text{subject to}\; q = \nabla u,\nonumber
\end{align}
where $q = (q_{1},q_{2})^{T}$ is the auxiliary term. An efficient solution to the problem shown in Eq. (\ref{eq:constrainedTV}) is using the augmented Lagrangian method defined as~\cite{Nocedal2006,Tai2009,Wu2010}
\begin{multline}\label{eq:ALM}
\underset{u,q}{\min}\;\underset{\mu}{\max}\; L_{TV}\left(u,q,\mu\right) = \frac{\lambda}{2}\int_{\Omega}\left(u - f\right)^{2}d\mathbf{x} + \int_{\Omega}\left|q\right| d\mathbf{x} \\ 
+ \int_{\Omega}\mu\cdot\left(q-\nabla u\right) d\mathbf{x} + \frac{r}{2}\int_{\Omega}\left|q-\nabla u\right|^{2} d\mathbf{x},
\end{multline}
where $\mu = (\mu_{1},\mu_{2})^{T}$ is the vector of Lagrange multipliers and $r$ is a positive constant. The iterative process to solve Eq. (\ref{eq:ALM}) is sketched in Algorithm \ref{alg:almTV}.~\cite{Nocedal2006,Tai2009,Wu2010}
\begin{algorithm}[ht!]
\DontPrintSemicolon
\KwData{$u^{0} = 0,\; q^{0} = 0,\; \mu^{0} = 0$}
$k = 0$\;
\While{stop criteria is not fulfilled}
  {
  \begin{align}
  \text{Solve}\qquad\left(u^{k+1},q^{k+1}\right) &\approx  \underset{u,q}{\min}\; L_{TV}(u,q,\mu^{k}; f)\label{eq:subproblema1}\\
  \text{then update}\qquad\mu^{k+1} &= \mu^{k} + r\left(q^{k+1} - \nabla u^{k+1}\right)\label{eq:subproblema2}
  \end{align}\;
  $k = k + 1$
  }\caption{Augmented Lagrangian method for the TV model}\label{alg:almTV}
\end{algorithm}

Eq. (\ref{eq:subproblema1}) is an unconstrained optimization problem which is difficult to solve because variables $u$ and $q$ are coupled. One alternative, proposed in references \citenum{Tai2009} and \citenum{Wu2010}, is to separate Eq. (\ref{eq:subproblema1}) in two subproblems defined as
\begin{equation}\label{eq:funcional1}
\underset{u}{\min}\quad \frac{\lambda}{2}\int_{\Omega}\left|u -f\right|^{2}d\mathbf{x} - \int_{\Omega}\mu\cdot\nabla u\,d\mathbf{x} + \frac{r}{2}\int_{\Omega}\left|q - \nabla u\right|^{2}d\mathbf{x},
\end{equation}
given the term $q$, and
\begin{equation}\label{eq:funcional2}
 \underset{q}{\min}\quad \int_{\Omega}\left|q\right|d\mathbf{x} + \int_{\Omega}\mu\cdot q\,d\mathbf{x} + \frac{r}{2}\int_{\Omega}\left|q - \nabla u\right|^{2}d\mathbf{x},
 \end{equation} 
given the term $u.$

The optimality condition of the problem shown in Eq.(\ref{eq:funcional1}) gives a linear equation and can be solved efficiently using the Fourier Transform. On the other side, the problem shown in Eq.(\ref{eq:funcional2}) has a closed-form solution known as \textit{soft-thresholding operator}~\cite{Donoho1995,Wang2008a,Caboussat2009} and is defined as~\cite{Tai2009,Wu2010}
\begin{align}\label{eq:softOp}
  q =
    \begin{cases}
      \frac{1}{r}\left(1 - \frac{1}{|w|}\right)\,w, &\text{if}\; |w|> 1\\
      0,  &\text{if}\;  |w| \leq 1
    \end{cases}
\end{align}
where $w = r\nabla u - \mu^{k}.$

\subsection{Augmented Lagrangian method for computing discontinuous phase maps based on TV model}
Following the idea described previously, we transform the problem shown in Eq. (\ref{eq:energiaMinimiza}) into a constrained one, and then solve it with the augmented Lagrangian method. An obvious advantage of using the augmented Lagrangian method is that the solution benefits from a fast solver and closed-form solutions. In our case, unfortunately some terms in Eq. (\ref{eq:energiaMinimiza}) are nonlinear, so a fast solver cannot be used in our solution; however the described closed-form solutions are susceptible to be used in our approach.

The proposed augmented Lagrangian method for Eq. (\ref{eq:energiaMinimiza}) is defined as
\begin{align}\label{eq:funcionalALM}
&\underset{\phi,a,b,\mathbf{q}}{\min}\;\underset{\mu}{\max}\; L\left(\phi_{\mathbf{x}}, b_{\mathbf{x}}, a_{\mathbf{x}},\mathbf{q}_{\phi},\mathbf{q}_{b},\mathbf{q}_{a},\mu_{\phi},\mu_{b},\mu_{a};\omega_{\mathbf{x}}\right) =\nonumber\\ 
&\frac{\lambda}{2}\int_{\Omega}\left(I_{\mathbf{x}} - g_{\mathbf{x}}\right)^{2}d\mathbf{x} + \int_{\Omega}\left|\mathbf{q}_{\phi}\right| d\mathbf{x} + \int_{\Omega}\left|\mathbf{q}_{b}\right| d\mathbf{x} + \int_{\Omega}\left|\mathbf{q}_{a}\right| d\mathbf{x}\nonumber\\
&+ \int_{\Omega}\mu_{\phi}\cdot\left(\mathbf{q}_{\phi}-\nabla\phi_{\mathbf{x}}\right) d\mathbf{x} + \int_{\Omega}\mu_{b}\cdot\left(\mathbf{q}_{b}-\nabla b_{\mathbf{x}}\right) d\mathbf{x}\nonumber\\ 
&+ \int_{\Omega}\mu_{a}\cdot\left(\mathbf{q}_{a}-\nabla a_{\mathbf{x}}\right) d\mathbf{x} + \frac{r}{2}\int_{\Omega}\left|\mathbf{q}_{\phi}-\nabla\phi_{\mathbf{x}}\right|^{2} d\mathbf{x}\nonumber\\ 
&+ \frac{r}{2}\int_{\Omega}\left|\mathbf{q}_{b}-\nabla b_{\mathbf{x}}\right|^{2} d\mathbf{x} + \frac{r}{2}\int_{\Omega}\left|\mathbf{q}_{a}-\nabla a_{\mathbf{x}}\right|^{2} d\mathbf{x},
\end{align}
where $r$ is a positive constant, $\mathbf{q}_{d} = (q_{1},q_{2})_{d}^{T}$ is the auxiliary term,  $\mu _{d}= (\mu_{1},\mu_{2})_{d}^{T}$ are Lagrange multipliers and $d$ is any variable representing $\phi,$ $b,$ or $a.$

As can be observed, the functional shown in Eq. (\ref{eq:funcionalALM}) has similar structure than the one shown in Eq. (\ref{eq:ALM}), so we follow the procedure described in the previous section to propose the solution of Eq. (\ref{eq:funcionalALM}). 

The minimization problem given in Eq. (\ref{eq:funcionalALM}) is separated in two subproblems. The first subproblem is related to the solution of the terms $\phi_{\mathbf{x}},$ $b_{\mathbf{x}},$  and $a_{\mathbf{x}}$, given the auxiliary terms $\mathbf{q}_{\phi},$ $\mathbf{q}_{b},$ and $\mathbf{q}_{a}.$ This subproblem is defined as
\begin{align}\label{eq:subprobPBA}
\underset{\phi}{\min}\; &\frac{\lambda}{2}\int_{\Omega}\left(I_{\mathbf{x}} - g_{\mathbf{x}}\right)^{2}d\mathbf{x} 
- \int_{\Omega}\mu_{\phi}\cdot\nabla\phi_{\mathbf{x}} d\mathbf{x}
+ \frac{r}{2}\int_{\Omega}\left|\mathbf{q}_{\phi}-\nabla\phi_{\mathbf{x}}\right|^{2} d\mathbf{x},\nonumber\\
\underset{b}{\min}\; &\frac{\lambda}{2}\int_{\Omega}\left(I_{\mathbf{x}} - g_{\mathbf{x}}\right)^{2}d\mathbf{x} 
- \int_{\Omega}\mu_{b}\cdot\nabla b_{\mathbf{x}} d\mathbf{x}
+ \frac{r}{2}\int_{\Omega}\left|\mathbf{q}_{b}-\nabla b_{\mathbf{x}}\right|^{2} d\mathbf{x},\\
\underset{a}{\min}\; &\frac{\lambda}{2}\int_{\Omega}\left(I_{\mathbf{x}} - g_{\mathbf{x}}\right)^{2}d\mathbf{x} 
- \int_{\Omega}\mu_{a}\cdot\nabla a_{\mathbf{x}} d\mathbf{x}
+ \frac{r}{2}\int_{\Omega}\left|\mathbf{q}_{a}-\nabla a_{\mathbf{x}}\right|^{2} d\mathbf{x}.\nonumber
\end{align}
The second subproblem is related to the solution of the auxiliary terms $\mathbf{q}_{\phi},$ $\mathbf{q}_{b},$ and $\mathbf{q}_{a},$ given the terms $\phi_{\mathbf{x}},$ $b_{\mathbf{x}},$  and $a_{\mathbf{x}}.$ This subproblem is defined as 
\begin{align}\label{eq:subprobQ}
\underset{\mathbf{q}_{\phi}}{\min}\; &\int_{\Omega}\left|\mathbf{q}_{\phi}\right|d\mathbf{x} + \int_{\Omega}\mu_{\phi}\cdot \mathbf{q}_{\phi}\,d\mathbf{x} + \frac{r}{2}\int_{\Omega}\left|\mathbf{q}_{\phi} - \nabla\phi_{\mathbf{x}}\right|^{2}d\mathbf{x},\nonumber\\
\underset{\mathbf{q}_{b}}{\min}\; &\int_{\Omega}\left|\mathbf{q}_{b}\right|d\mathbf{x} + \int_{\Omega}\mu_{b}\cdot \mathbf{q}_{b}\,d\mathbf{x} + \frac{r}{2}\int_{\Omega}\left|\mathbf{q}_{b} - \nabla b_{\mathbf{x}}\right|^{2}d\mathbf{x},\\
\underset{\mathbf{q}_{a}}{\min}\; &\int_{\Omega}\left|\mathbf{q}_{a}\right|d\mathbf{x} + \int_{\Omega}\mu_{a}\cdot \mathbf{q}_{a}\,d\mathbf{x} + \frac{r}{2}\int_{\Omega}\left|\mathbf{q}_{a} - \nabla a_{\mathbf{x}}\right|^{2}d\mathbf{x}.\nonumber
\end{align}

In the case of the functionals shown in Eq. (\ref{eq:subprobPBA}), the solution can be stated as follows: 

The first-order optimality conditions of Eq. (\ref{eq:subprobPBA}) are given by
\begin{align}
-&\lambda\Big(a_{\mathbf{x}} + b_{\mathbf{x}}\cos\psi_{\mathbf{x}} - g_{\mathbf{x}}\Big)\Big(b_{\mathbf{x}}\sin\psi_{\mathbf{x}}\Big) + \nabla\cdot\mu_{\phi} + r\nabla\cdot\Big(\mathbf{q}_{\phi}-\nabla\phi_{\mathbf{x}}\Big) = 0,\nonumber\\
&\lambda\Big(a_{\mathbf{x}} + b_{\mathbf{x}}\cos\psi_{\mathbf{x}} - g_{\mathbf{x}}\Big)\Big(\cos\psi_{\mathbf{x}}\Big) + \nabla\cdot\mu_{b} + r\nabla\cdot\Big(\mathbf{q}_{b}-\nabla b_{\mathbf{x}}\Big) = 0,\nonumber\\
&\lambda\Big(a_{\mathbf{x}} + b_{\mathbf{x}}\cos\psi_{\mathbf{x}} - g_{\mathbf{x}}\Big) + \nabla\cdot\mu_{a} 
+ r\nabla\cdot\Big(\mathbf{q}_{a} - \nabla a_{\mathbf{x}}\Big) = 0.
\end{align}
with boundary conditions
\begin{gather*}
\mu_{\phi}\cdot\mathbf{n} = 0,\quad \Big(\mathbf{q}_{\phi}-\nabla\phi_{\mathbf{x}}\Big)\cdot\mathbf{n} = 0,\\
\mu_{b}\cdot\mathbf{n} = 0,\quad \Big(\mathbf{q}_{b}-\nabla b_{\mathbf{x}}\Big)\cdot\mathbf{n} = 0,\\
\mu_{a}\cdot\mathbf{n} = 0,\quad \Big(\mathbf{q}_{a}-\nabla a_{\mathbf{x}}\Big)\cdot\mathbf{n} = 0.
\end{gather*}

The above equations have the same form than Eq. (\ref{eq:gradienteE1}), so it is possible to express them as the fixed point iteration shown in Eq. (\ref{eq:fixedpoint}):
\begin{align}\label{eq:fixedpointALM}
  \Big(\lambda \left(b^{k}_{\mathbf{x}}\right)^{2}\sin^{2}\psi^{k}_{\mathbf{x}} - r\nabla\cdot\nabla\Big)\phi^{k+1}_{\mathbf{x}} &= \lambda a^{k}_{\mathbf{x}}b^{k}_{\mathbf{x}}\sin\psi^{k}_{\mathbf{x}}\nonumber\\
  + \lambda\left(b^{k}_{\mathbf{x}}\right)^{2}\cos\psi^{k}_{\mathbf{x}}\sin\psi^{k}_{\mathbf{x}} &+ \lambda\phi^{k}_{\mathbf{x}}\left(b^{k}_{\mathbf{x}}\right)^{2}\sin^{2}\psi^{k}_{\mathbf{x}}\nonumber\\
  - \lambda g_{\mathbf{x}}b^{k}_{\mathbf{x}}\sin\psi^{k}_{\mathbf{x}} &- \nabla\cdot\mu_{\phi}^{k} - r\nabla\cdot\mathbf{q}_{\phi}^{k},\nonumber\\
  \Big(\lambda\cos^{2}\psi^{k}_{\mathbf{x}} - r\nabla\cdot\nabla\Big) b^{k+1}_{\mathbf{x}} &= -\lambda\Big(a^{k}_{\mathbf{x}}\cos\psi^{k}_{\mathbf{x}} - g_{\mathbf{x}}\cos\psi^{k}_{\mathbf{x}}\Big)\nonumber\\ 
  - \nabla\cdot\mu_{b}^{k} &- r\nabla\cdot\mathbf{q}_{b}^{k}\\
  \Big(\lambda - r\nabla\cdot\nabla\Big)a^{k+1}_{\mathbf{x}} = -\lambda\Big(b^{k}_{\mathbf{x}}\cos\psi^{k}_{\mathbf{x}} &- g_{\mathbf{x}}\Big) - \nabla\cdot\mu_{a}^{k} - r\nabla\cdot\mathbf{q}_{a}^{k}.\nonumber
\end{align}

On the other hand, in the case of the functionals shown in Eq. (\ref{eq:subprobQ}), we found these to have the same structure that in Eq.(\ref{eq:funcional2}), so we can use the \textit{soft-thresholding operator}~\cite{Donoho1995,Wang2008a,Caboussat2009} to solve them:
\begin{align}\label{eq:thresholdingALM}
  \mathbf{q}_{d}^{k+1} =
    \begin{cases}
      \frac{1}{r}\left(1 - \frac{1}{|w_{d}|}\right)w_{d}, &\text{if}\; |w_{d}|> 1\\
      0,  &\text{if}\;  |w_{d}| \leq 1
    \end{cases}    
\end{align}
where $w_{d} = r\nabla d_{\mathbf{x}}^{k+1} - \mu_{d}^{k},$ and $d$ is any variable representing $\phi,$ $b,$ or $a.$

Finally, the update of the Lagrange multipliers is carried out in the same way it was done in Eq. (\ref{eq:subproblema2}). The iterative procedure to solve Eq. (\ref{eq:funcionalALM}) is given in Algorithm \ref{alg:algoritmo}. 
\begin{algorithm}[ht!]
\DontPrintSemicolon
\KwData{$\phi^{0}_{\mathbf{x}} = 0,\;b^{0}_{\mathbf{x}} = 0,\;a^{0}_{\mathbf{x}} = 0,\;\mathbf{q}^{0}_{\phi} = 0,\;\mathbf{q}^{0}_{b} = 0,\;\mathbf{q}^{0}_{a} = 0,$\\\begin{center} $\mu^{0}_{\phi} = 0,\;\mu^{0}_{b} = 0,\;\mu^{0}_{a} = 0$\end{center}}
$k = 0$\;
\While{stop criteria is not fulfilled}
  {\;
    Compute $\phi^{k+1}_{\mathbf{x}},$ $b^{k+1}_{\mathbf{x}},$ and $a^{k+1}_{\mathbf{x}}$ using the fixed point iteration shown in Eq. (\ref{eq:fixedpointALM}).\;\;
    Compute $\mathbf{q}_{\phi}^{k+1},$ $\mathbf{q}_{b}^{k+1},$ and $\mathbf{q}_{a}^{k+1}$ using Eq. (\ref{eq:thresholdingALM}).\;\;
    Update
    \begin{align}
      \mu_{\phi}^{k+1} = \mu_{\phi}^{k} + r\left(\mathbf{q}_{\phi}^{k+1} - \nabla \phi_{\mathbf{x}}^{k+1}\right)\nonumber\\%\;\;\text{where}\; d = \phi, b, a
      \mu_{b}^{k+1} = \mu_{b}^{k} + r\left(\mathbf{q}_{b}^{k+1} - \nabla b_{\mathbf{x}}^{k+1}\right)\\
      \mu_{a}^{k+1} = \mu_{a}^{k} + r\left(\mathbf{q}_{a}^{k+1} - \nabla a_{\mathbf{x}}^{k+1}\right)\nonumber
    \end{align}
    $k = k + 1$\;
  }\caption{Augmented Lagrangian method for Eq. (\ref{eq:funcionalALM})}\label{alg:algoritmo}
\end{algorithm}

\section{Numerical Experiments}
To illustrate the performance of the augmented Lagrangian method (ALM), we carried out some numerical experiments using a Intel Core i7 @ 2.40 GHz laptop with Debian GNU/Linux 9 (64-bit) and 16 GB of memory. In these experiments, we compare our proposed method with the fixed point method (FP) shown in Eq. (\ref{eq:fixedpoint}). Both methods were implemented in C/C++. In our experiments we used as stopping criteria the following condition \[\frac{\|d^{k} -d^{k-1}\|}{\|d^{k-1}\|} \leq \epsilon\] where $\epsilon = 10^{-5}$ and $d$ is any variable representing $\phi,$ $b,$ or $a.$ For the augmented Lagrangian method we use $r = 11.5,$ and for the fixed point method we use $\beta = 10^{-3},$ which is the constant to avoid division by zero.

For simplicity, we selected  the regularization parameter $\lambda$ manually.  However, well known methods can be used to obtain the best parameter for this task, such as those described in section 5.6 of reference \citenum{Bertero1998}. In addition, we use a normalized error $Q$ to compare the phase-map estimation; this error is defined as~\cite{Perlin2016} \[Q\left(\mu,\nu\right) = \frac{\| \mu - \nu\|_{2}}{\| \mu\|_{2} + \|\nu\|_{2}},\] where $\mu$ and $\nu$ are the signals to be compared. The normalized error values vary between zero (for perfect agreement) and one (for total disagreement).

\subsection{Phase demodulation using synthetic fringe pattern}
\begin{figure}[ht!]
  \centering
  \includegraphics[width=\textwidth]{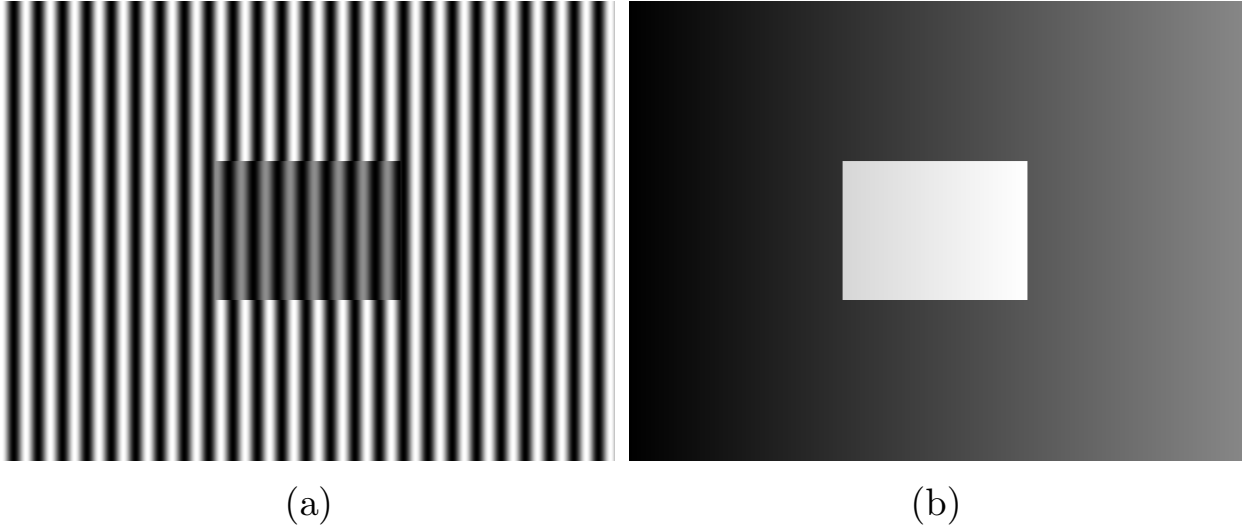}
  \caption{(a) Synthetic fringe pattern. (b) Synthetic phase term $\phi_{\mathbf{x}}.$}
  \label{fig:f1synthetic}
\end{figure}

Here we present two experiments using a synthetic fringe pattern of size 640 x 480 pixels, generated in similar way to that described in Refs. \citenum{Galvan2006} and \citenum{Legarda-Saenz2014}. Figure \ref{fig:f1synthetic} shows the synthetic fringe pattern and the synthetic phase term $\phi_{\mathbf{x}}.$ 

\begin{figure}[ht!]
  \centering
  \includegraphics[width=\textwidth]{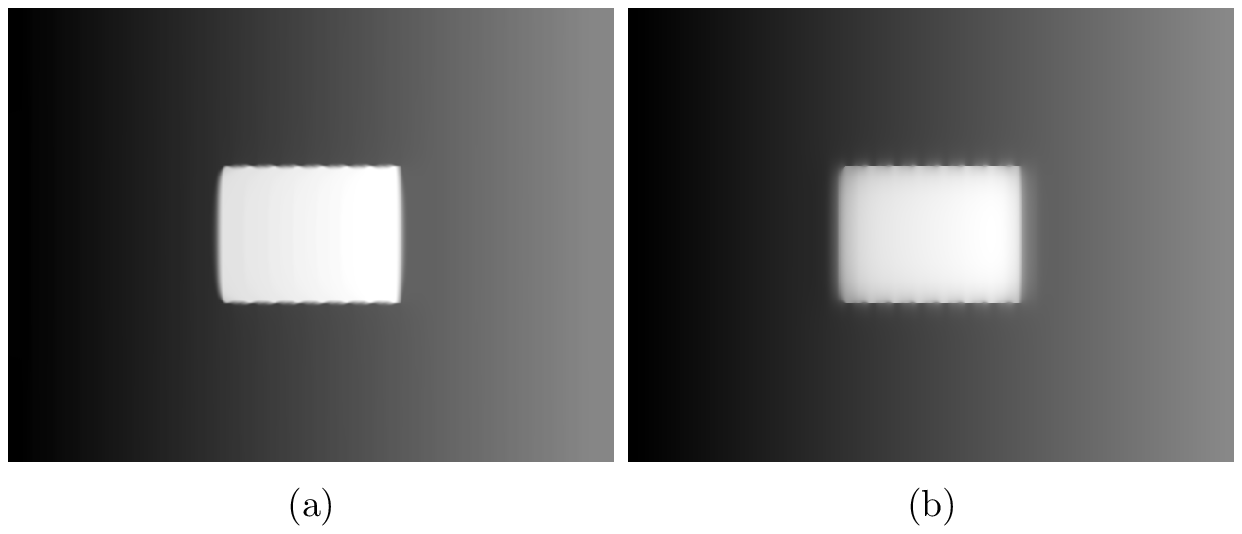}
  \caption{Estimated phase terms using (a) augmented Lagrangian method, Eq. (\ref{eq:funcionalALM}), and (b) fixed point method, Eq. (\ref{eq:fixedpoint}).}
  \label{fig:f2estimation}
\end{figure}

\begin{figure}[ht!]
  \centering
  \includegraphics[width=\textwidth]{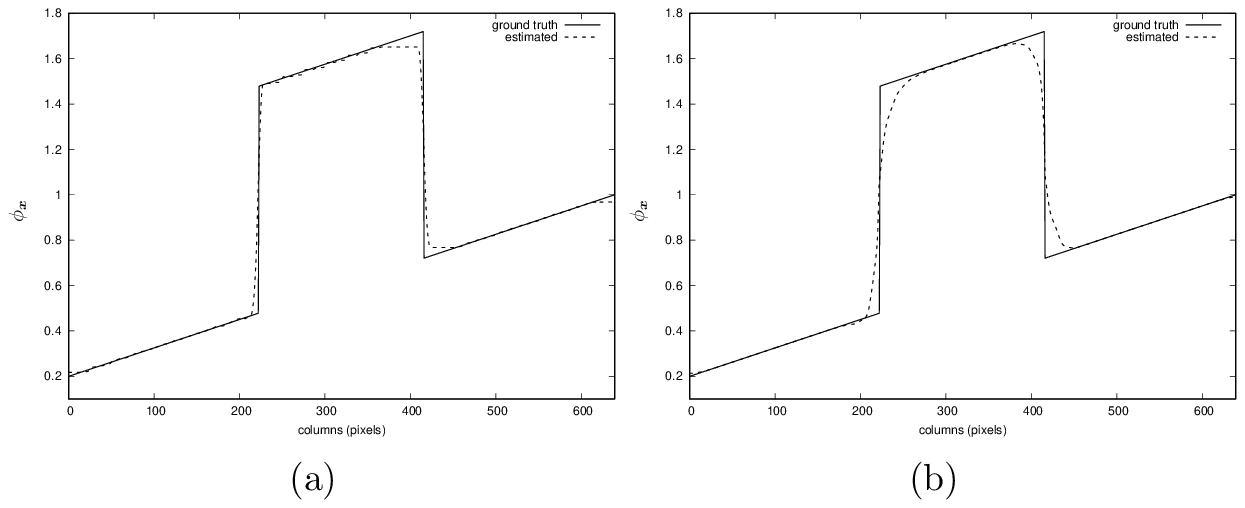}
  \caption{Middle row of the estimated phase terms using (a) augmented Lagrangian method, Eq. (\ref{eq:funcionalALM}), and (b) fixed point method, Eq. (\ref{eq:fixedpoint}).}
  \label{fig:f3estimation}
\end{figure}

The first experiment was the demodulation of this fringe pattern using the augmented Lagrangian method and the fixed point method, both with $\lambda = 10.$ The resultant phase demodulations are shown in Figure \ref{fig:f2estimation}. In Figure \ref{fig:f3estimation} we show the middle row of the estimated phase term $\phi_{\mathbf{x}}.$ The normalized error of the augmented Lagrangian method was $Q = 0.0243$ and the time employed to obtain the solution was 141 seconds using 655 iterations. On the other hand, the normalized error of the fixed point method was $Q = 0.0343$ and the time employed to obtain the solution was 340 seconds using 841 iterations. 

In Figure \ref{fig:f2estimation}, it can be seen that both methods successfully demodulate the discontinuity found in the synthetic fringe pattern. However, when analysing Figure \ref{fig:f3estimation}, we found that the demodulation of the proposed augmented Lagrangian method is more precise than that of the fixed point method. Moreover, the time and number of iterations employed to obtain the solution are better than the fixed point method. This is due to the influence of the term $\beta$ on the performance of the fixed point method: with a larger value of this term, the speed of convergence of the method is better but the accuracy of the solution gets worse. This is not desirable in fringe analysis techniques.

\begin{figure}[ht!]
  \centering
  \includegraphics[width=\textwidth]{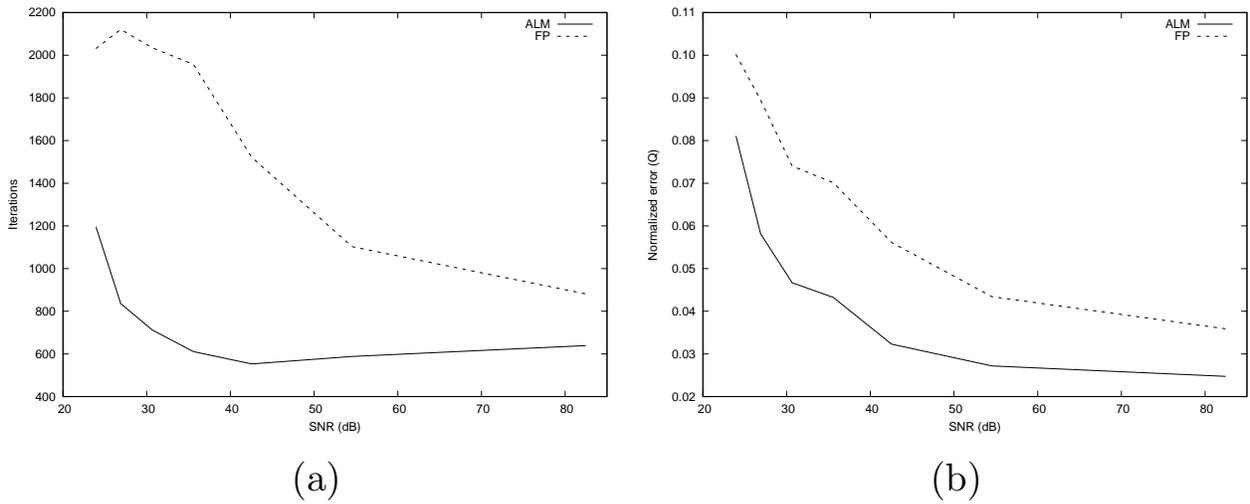}
  \caption{Performance of the phase demodulation methods with different noise levels: (a) iterations employed, and (b) the normalized error $Q.$}
  \label{fig:f4estimation}
\end{figure}

The second experiment was the demodulation of the fringe pattern shown in Figure \ref{fig:f2estimation} with different levels of noise. In this experiment we used $\lambda = 6$ for both methods. The resultant performance of both methods is shown in Figure \ref{fig:f4estimation}. As can be observed, in this experiment we found the same differences mentioned previously: the proposed augmented Lagrangian method demodulates the fringe pattern faster and more accurately than the fixed point method, even with noisy fringe patterns.

\subsection{Phase demodulation using experimental fringe patterns}
\begin{figure}[ht!]
  \centering
  \includegraphics[width=\textwidth]{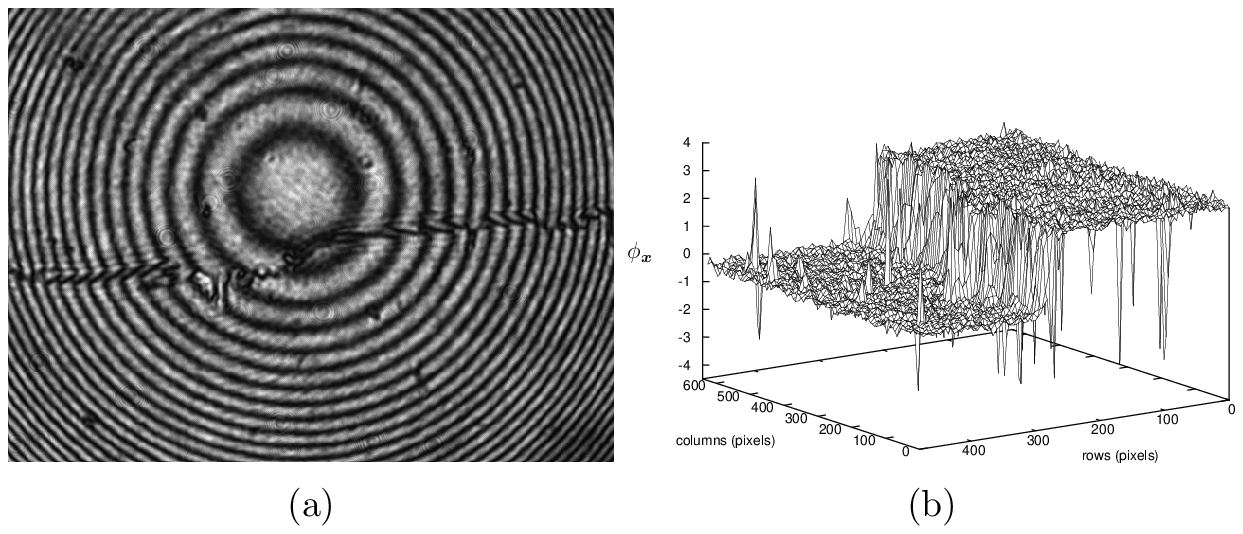}
  \caption{(a) Experimental fringe pattern. (b) Estimated phase term $\phi_{\mathbf{x}}$ using Schwider-Hariharan (4+1) algorithm.~\cite{Kreis1996,Servin2014}}
  \label{fig:f5reales}
\end{figure}

Here we present the phase demodulation of a fringe pattern obtained from a holographic interferometry experiment~\cite{Kreis1996}, which consisted of the height measurement of a micro-thin film. The fringe pattern obtained from this experiment, with 640 x 480 pixels, is shown in Figure \ref{fig:f5reales}, panel (a). Figure \ref{fig:f5reales}, panel (b) shows a phase term $\phi_{\mathbf{x}}$ estimation of this experimental fringe pattern using Schwider-Hariharan (4+1) algorithm.~\cite{Kreis1996,Servin2014} This estimation was used as reference in this experiment.

\begin{figure}[ht!]
  \centering
  \includegraphics[width=\textwidth]{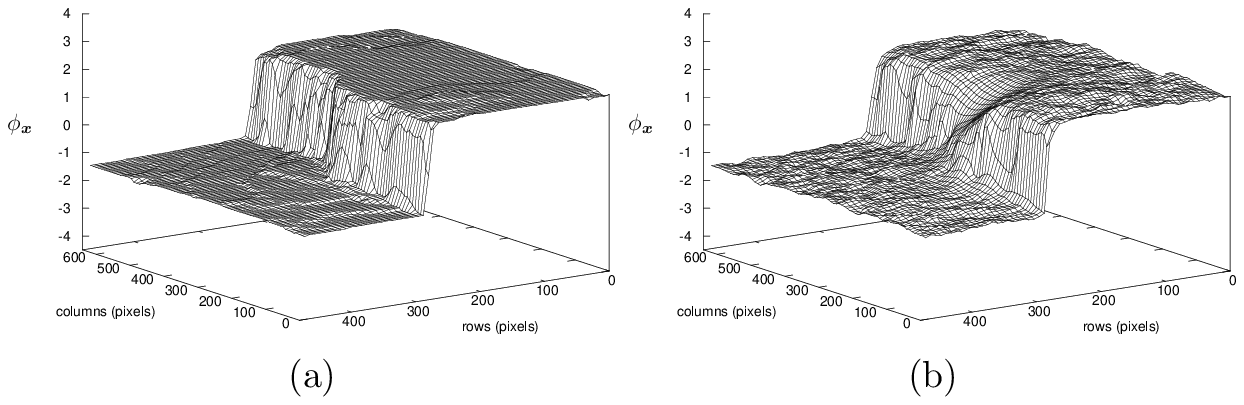}
  \caption{Estimated phase terms using (a) augmented Lagrangian method, Eq. (\ref{eq:funcionalALM}), and (b) fixed point method, Eq. (\ref{eq:fixedpoint}).}
  \label{fig:f6estimareales}
\end{figure}

\begin{figure}[ht!]
  \centering
  \includegraphics[width=\textwidth]{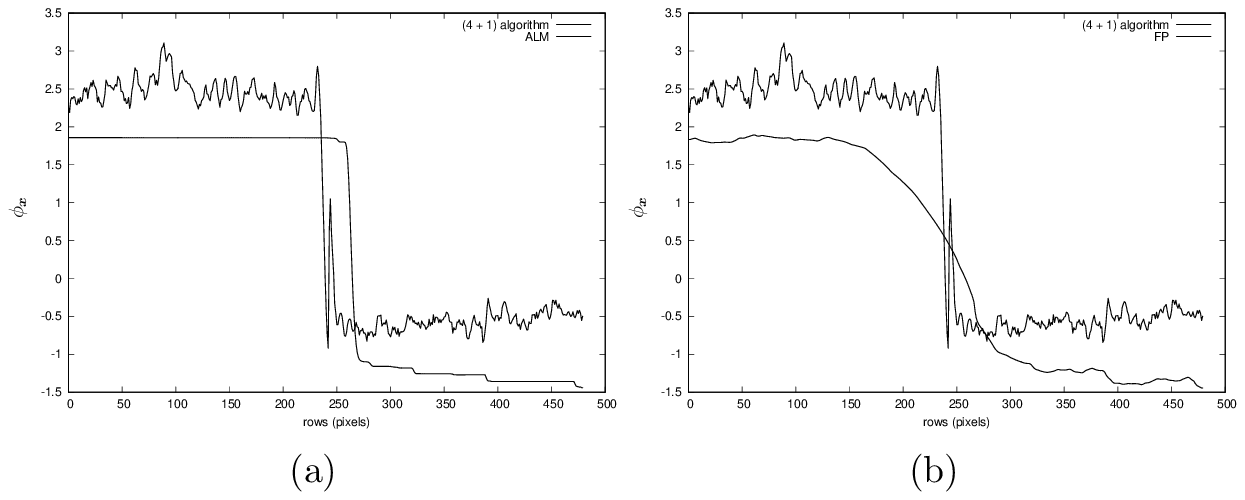}
  \caption{Middle column of the estimated phase terms using (a) augmented Lagrangian method, Eq. (\ref{eq:funcionalALM}), and (b) fixed point method, Eq. (\ref{eq:fixedpoint}).}
  \label{fig:f7estimareales}
\end{figure}

The demodulations of this experimental fringe pattern using the augmented Lagrangian method and the fixed point method are shown in Figure \ref{fig:f6estimareales}, both with $\lambda = 10.$ In Figure \ref{fig:f7estimareales} we shows the middle column of the estimated phase term $\phi_{\mathbf{x}}.$ The time employed by the augmented Lagrangian method to obtain the solution was 430 seconds using 1995 iterations. On the other hand, the fixed point method took 3288 seconds using 7445 iterations to obtain the solution. 

As can be observed in Figures \ref{fig:f6estimareales} and \ref{fig:f7estimareales}, both methods are able to demodulate the discontinuity found in the experimental fringe pattern. These estimations can be seen as the filtered version of the one shown in Figure \ref{fig:f5reales}, panel (b). One relevant aspect is that both methods preserve the dynamic range of the phase term $\phi_{\mathbf{x}}.$ On the other hand, there are marked differences in the quality of the estimation and in the numerical performance to obtain it: the proposed augmented Lagrangian method delivers a piecewise constant surface while the fixed point method gives a smoother one. Additionally, our proposal is several times faster than the fixed point method.

\section{Conclusions}
In this paper we present a reformulation of the method presented in reference \citenum{Legarda-Saenz2014} as a constrained minimization problem and we solve it using the augmented Lagrangian method. As can be seen in the numerical experiments, the proposed method is able to accurately demodulate a single fringe pattern with discontinuities. The numerical solution of Eq. (\ref{eq:funcionalALM}) results in a simple algorithm, which is faster and preserves the dynamic range of the phase term $\phi_{\mathbf{x}}.$ An extra advantage of the proposed method is its feasibility to be implemented on dedicated hardware to obtain real-time processing. This will be one aim of our future research.

%%%%%%%%%%%%%%%%%%%%%%%%%%%%%%%%
%           Referencias
%%%%%%%%%%%%%%%%%%%%%%%%%%%%%%%%
\newpage
\bibliographystyle{unsrt}

\begin{thebibliography}{10}

\bibitem{Servin2014}
Manuel Servin, Juan~Antonio Quiroga, and Moises Padilla.
\newblock {\em {Fringe Pattern Analysis for Optical Metrology: Theory,
  Algorithms, and Applications}}.
\newblock Wiley-VCH, Weinheim, 2014.

\bibitem{Marroquin1999}
J.~L. Marroquin, M.~Rivera, S.~Botello, R.~Rodriguez-Vera, and M.~Servin.
\newblock {Regularization methods for processing fringe-pattern images}.
\newblock {\em Applied Optics}, 38(5):788--794, 1999.

\bibitem{Legarda-Saenz2002a}
Ricardo Legarda-Saenz, Wolfgang Osten, and Werner~P. Juptner.
\newblock {Improvement of the regularized phase tracking technique for the
  processing of nonnormalized fringe patterns}.
\newblock {\em Applied Optics}, 41(26):5519--5526, 2002.

\bibitem{Villa2000}
Jesus Villa, Juan~Antonio Quiroga, and Manuel Servin.
\newblock {Improved regularized phase-tracking technique for the processing of
  squared-grating deflectograms}.
\newblock {\em Applied Optics}, 39(4):502--508, 2000.

\bibitem{Rivera2005}
Mariano Rivera.
\newblock {Robust phase demodulation of interferograms with open or closed
  fringes}.
\newblock {\em Journal of the Optical Society of America A}, 22(6):1170--1175,
  2005.

\bibitem{Galvan2006}
Carlos Galvan and Mariano Rivera.
\newblock {Second-order robust regularization cost function for detecting and
  reconstructing phase discontinuities}.
\newblock {\em Applied Optics}, 45(2):353--359, 2006.

\bibitem{Rudin1992}
Leonid~I. Rudin, Stanley Osher, and Emad Fatemi.
\newblock {Nonlinear total variation based noise removal algorithms}.
\newblock {\em Physica D: Nonlinear Phenomena}, 60(1-4):259--268, 1992.

\bibitem{Legarda-Saenz2014}
Ricardo Legarda-Saenz, Carlos Brito-Loeza, and Arturo Espinosa-Romero.
\newblock {Total variation regularization cost function for demodulating phase
  discontinuities}.
\newblock {\em Applied Optics}, 53(11):2297, 2014.

\bibitem{Brito-Loeza2019}
Carlos Brito‐Loeza, Ricardo Legarda‐Saenz, and Anabel Martin‐Gonzalez.
\newblock {A fast algorithm for a total variation based phase demodulation
  model}.
\newblock {\em Numerical Methods for Partial Differential Equations}, n/a(n/a),
  nov 2019.

\bibitem{Vogel1996}
C.~R. Vogel and M.~E. Oman.
\newblock {Iterative methods for total variation denoising}.
\newblock {\em SIAM Journal on Scientific Computing}, 17(1):227--238, 1996.

\bibitem{Vogel2002}
Curtis~R. Vogel.
\newblock {\em {Computational Methods for Inverse Problems}}.
\newblock SIAM, 2002.

\bibitem{Chan1999}
Tony~F. Chan, Gene~H. Golub, and Pep Mulet.
\newblock {A Nonlinear Primal-Dual Method for Total Variation-Based Image
  Restoration}.
\newblock {\em SIAM Journal on Scientific Computing}, 20(6):1964--1977, jan
  1999.

\bibitem{Chambolle2004}
Antonin Chambolle.
\newblock {An algorithm for total variation minimization and applications}.
\newblock {\em Journal of Mathematical Imaging and Vision}, 20(1):89--97, 2004.

\bibitem{Getreuer2012}
Pascal Getreuer.
\newblock {Rudin-Osher-Fatemi Total Variation Denoising using Split Bregman}.
\newblock {\em Image Processing On Line}, 2012:1--20, may 2012.

\bibitem{Hestenes1969}
Magnus~R. Hestenes.
\newblock {Multiplier and gradient methods}.
\newblock {\em Journal of Optimization Theory and Applications}, 4(5):303--320,
  nov 1969.

\bibitem{Rockafellar1973}
R.~Tyrrell Rockafellar.
\newblock {A dual approach to solving nonlinear programming problems by
  unconstrained optimization}.
\newblock {\em Mathematical Programming}, 5(1):354--373, dec 1973.

\bibitem{Bertsekas1996}
Dimitri~P. Bertsekas.
\newblock {\em {Constrained Optimization and Lagrange Multiplier Methods}}.
\newblock Athena Scientific, Belmont, MA, 1996.

\bibitem{Nocedal2006}
Jorge Nocedal and Stephen Wright.
\newblock {\em {Numerical Optimization}}.
\newblock Springer, New York, second edition, 2006.

\bibitem{Tai2009}
Xue-Cheng Tai and Chunlin Wu.
\newblock {Augmented Lagrangian method, dual methods and split Bregman
  iteration for ROF model}.
\newblock In {\em Second International Conference, SSVM 2009}, volume LNCS
  5567, pages 502--513, 2009.

\bibitem{Wu2010}
Chunlin Wu and Xue-Cheng Tai.
\newblock {Augmented Lagrangian method, dual methods, and split Bregman
  iteration for ROF, vectorial TV, and high order models}.
\newblock {\em SIAM Journal on Imaging Sciences}, 3(3):300--339, 2010.

\bibitem{Donoho1995}
D.L. Donoho.
\newblock {De-noising by soft-thresholding}.
\newblock {\em IEEE Transactions on Information Theory}, 41(3):613--627, may
  1995.

\bibitem{Wang2008a}
Yilun Wang, Junfeng Yang, Wotao Yin, and Yin Zhang.
\newblock {A New Alternating Minimization Algorithm for Total Variation Image
  Reconstruction}.
\newblock {\em SIAM Journal on Imaging Sciences}, 1(3):248--272, jan 2008.

\bibitem{Caboussat2009}
A.~Caboussat, R.~Glowinski, and V.~Pons.
\newblock {An augmented Lagrangian approach to the numerical solution of a
  non-smooth eigenvalue problem}.
\newblock {\em Journal of Numerical Mathematics}, 17(1):3--26, jan 2009.

\bibitem{Bertero1998}
M.~Bertero and P.~Boccacci.
\newblock {\em {Introduction to Inverse Problems in Imaging}}.
\newblock Institute of Physics Publishing, Bristol, 1998.

\bibitem{Perlin2016}
Marc Perlin and Miguel~D. Bustamante.
\newblock {A robust quantitative comparison criterion of two signals based on
  the Sobolev norm of their difference}.
\newblock {\em Journal of Engineering Mathematics}, 101(1):115--124, dec 2016.

\bibitem{Kreis1996}
Thomas Kreis.
\newblock {\em {Holographic Interferometry: Principles and Methods}}.
\newblock Wiley-VCH, Berlin, 1996.

\end{thebibliography}

%%%%%%%%%%%%%%%%%%%%%%%%%%%%%%%%
%           Figuras
%%%%%%%%%%%%%%%%%%%%%%%%%%%%%%%%
\newpage
\section*{Figures}

\end{document}